\documentclass[11pt]{amsart}

\usepackage{array}
\newcolumntype{C}{>{$}c<{$}}
\newcolumntype{L}{>{$}l<{$}}
\newcolumntype{R}{>{$}r<{$}}

\newtheorem{theorem}{Theorem}

\theoremstyle{definition}

\newcommand{\C}{\mathbb{C}}

\newcommand{\Q}{\mathbb{Q}}

\DeclareMathOperator{\Exp}{\mathsf{Exp}}
\DeclareMathOperator{\Log}{\mathsf{Log}}

\renewcommand{\SS}{\mathbb{S}}
\DeclareMathOperator{\Serre}{\mathsf{e}}
\newcommand{\Mbar}{\overline{\mathcal{M}}}
\newcommand{\CM}{\mathcal{M}}
\renewcommand{\]}{{]\!]}}
\renewcommand{\[}{{[\![}}

\renewcommand{\o}{\otimes}
\newcommand{\p}{\partial}

\begin{document}

\title{The Hodge polynomial of $\Mbar_{3,1}$}

\author{E. Getzler}

\address{Northwestern University, 2033 Sheridan Rd., Evanston, IL
60208-2730, USA}

\email{getzler@math.nwu.edu}

\author{E. Looijenga}

\address{Utrecht University, P.O. Box 80010, 3508 TA Utrecht, The
Netherlands}

\email{looijeng@math.uu.nl}

\thanks{The authors thank the Scuola Normale, Pisa, where this paper was
written. The work of the first author is partially funded by the NSF under
grant DMS-9704320}

\begin{abstract}
We give formulas for the Dolbeault numbers of $\Mbar_{3,1}$, using the
first author's calculations of the weights of the cohomology of $\CM_{2,2}$
and $\CM_{2,3}$ and the second author's calculation of the weights of the
cohomology of $\CM_{3,1}$.
\end{abstract}

\maketitle

\renewcommand{\baselinestretch}{1.05}

Recall that the Serre polynomial of a quasi-projective variety $X$, or more
generally, of a Deligne-Mumford stack whose coarse moduli space is a
quasi-projective variety, is the polynomial
$$
\Serre(X) = \sum_{i,j} u^i v^j \sum_{k=0}^\infty (-1)^k
h^{i,j}(H^k_c(X,\C)) .
$$
If a finite group $\Gamma$ acts on $X$, the equivariant Serre polynomial
$\Serre_\Gamma(X)$ is defined in the same way, and takes values in
$R(\Gamma)$. For background to all of this, see \cite{config}.

If $X$ is smooth and projective, or more generally, if $X$ is a
smooth Deligne-Mumford stack whose coarse moduli space is projective, then
the Serre polynomial of $X$ equals its Hodge polynomial:
$$
\Serre(X) = \sum_{i,j} (-u)^i (-v)^j h^{i,j}(X) .
$$
In particular, Poincar\'e duality implies the functional equation
$$
\Serre(X)(u,v) = q^{\dim(X)} \Serre(X)(u^{-1},v^{-1}) ,
$$
where throughout this paper, $q$ denotes $uv$. All of this is also true
equivariantly.

If $2g-2+n>0$, let $\CM_{g,n}$ be the moduli stack of smooth projective
curves of genus $g$ together with $n$ distinct marked points, and let
$\Mbar_{g,n}$ be the moduli stack of stable curves of (arithmetic) genus
$g$ with $n$ distinct marked smooth points \cite{Knudsen}; these stacks are
smooth, of dimension $3g-3+n$, and the coarse moduli space of $\Mbar_{g,n}$
is projective.  Furthermore, there is a flat map
$\Mbar_{g,n+1}\to\Mbar_{g,n}$ which represents $\Mbar_{g,n+1}$ as the
universal curve over $\Mbar_{g,n}$. In particular, $\Mbar_{g,1}$ is the
universal stable curve of genus $g$. In this paper, we calculate
$\Serre(\Mbar_{3,1})$. In particular, we show that $h^{i,j}(\Mbar_{3,1})=0$
unless $i=j$.

There is a standard identification of $R(\SS_n)$ with the space of
symmetric functions (in an infinite number of variables $\{x_i\mid
i\ge0\}$) of degree $n$, which we denote by $\Lambda_n$. If $X$ is a
quasi-projective variety with action of $\SS_n$, let $\Serre_n(X)$ be its
Serre polynomial, thought of as an element of $\Lambda_n[u,v]$.

There is a formula permitting the calculation of the $\SS_n$-equivariant
Hodge polynomial of $\Mbar_{g,n}$ in terms of the $\SS_m$-equivariant Serre
polynomials of the moduli stacks $\CM_{h,m}$, where $h\le g$ and $2h+m\le
2g+n$.

Let $\hat\Lambda$ be the completion of the ring of symmetric functions:
$$
\hat\Lambda = \prod_{n=0}^\infty \Lambda_n .
$$
We may identify $\hat\Lambda$ with the algebra of power series in the
variables
$$
h_n = \sum_{i_1\le\dots\le i_n} x_{i_1}\dots x_{i_n} .
$$
Introduce an auxilliary variable $\lambda$, and let
$\hat\Lambda\[\lambda\]$ be the ring of power series in $\lambda$ with
coefficients in the ring of symmetric functions. We may identify
$\hat\Lambda\o\Q$ with the algebra of power series in the variables
$$
p_n = \sum_{i=0}^\infty x_i^n .
$$
The ring $\hat\Lambda\[\lambda\]$ is a (special) lambda-ring: the Adams
operations on $\hat\Lambda\[\lambda\]$ are determined by the formulas
$\psi_k(p_n)=p_{kn}$ and $\psi_k(\lambda)=\lambda^k$. Let $\Exp$ be the
operation
$$
\Exp(f) = \sum_{k=0}^\infty \sigma_k(f) = \exp\left( \sum_{k=0}^\infty
\frac{\psi_k(f)}{k} \right) : \lambda \hat\Lambda\[\lambda\] \to 1 +
\lambda \hat\Lambda\[\lambda\] ,
$$
and let $\Log$ be its inverse, given by the formula
$$
\Log(f) = \sum_{k=0}^\infty \frac{\mu(k)}{k} \psi_k(f) : 1 + \lambda
\hat\Lambda\[\lambda\] \to \lambda \hat\Lambda\[\lambda\] .
$$
Let $\Delta$ be the differential operator on $\hat\Lambda$ (and by
extension of coefficients, on $\hat\Lambda\[\lambda\]$) give by the formula
$$
\Delta = \sum_{k=1}^\infty \lambda^{2k} \left( \frac{k}{2} \frac{\p^2}{\p
p_k^2} + \frac{\p}{\p p_{2k}} \right) .
$$

We may define two elements of $\hat\Lambda\[\lambda\]$, by the formulas
\begin{align*}
\Psi &= \sum_{2g-2+n>0,n>0} \lambda^{2g-2+n} \Serre_{\SS_n}(\Mbar_{g,n})
\quad \text{and} \\
\Phi &= \sum_{2g-2+n>0,n>0} \lambda^{2g-2+n} \Serre_{\SS_n}(\Mbar_{g,n}) .
\end{align*}
The following is Theorem (8.13) of \cite{modular}.
\begin{theorem} \label{Feynman}
$\Psi = \Log\bigl( \exp(\Delta) \Exp( \Phi ) \bigr)$
\end{theorem}

We see that to use this formula to calculate $\Serre(\Mbar_{3,1})$, we need
the Serre polynomial of $\CM_{3,1}$, together with the following
equivariant Serre polynomials:
$$\begin{tabular}{|C|L|} \hline
\CM_{g,n} & \Serre_{\SS_n}(\CM_{g,n}) \\ \hline
\CM_{0,3} & s_3 \\
\CM_{0,4} & qs_4-s_{2^2} \\
\CM_{0,5} & q^2s_5-qs_{32}+s_{31^2} \\
\CM_{0,6} & q^3s_6-q^2s_{42}+q(s_{41^2}+s_{321}) -
(s_{41^2}+s_{3^2}+s_{2^21^2}) \\
\CM_{0,7} & q^4s_7 - q^3s_{52} + q^2(s_{3^21}+s_{51^2}+s_{421}) \\
& \quad {} -
q(s_{51^2}+s_{43}+s_{421}+s_{2^31}+s_{41^3}+s_{3^21}+s_{321^2}) \\
& \quad {} + s_{321^2}+s_{52}+s_{31^4}+s_{32^2}+s_{421} \\[2pt] \hline
\CM_{1,1} & qs_1 \\
\CM_{1,2} & q^2s_2 \\
\CM_{1,3} & q^3s_3 - s_{1^3} \\
\CM_{1,4} & q^4s_4-q^2s_4-qs_{21^2}+s_{31} \\
\CM_{1,5} & q^5s_5-q^3(s_5+s_{41})+q^2(s_{32}-s_{31^2}) \\
& \quad {} + q(s_{41}+s_{32}+s_{31^2}) - (s_5+s_{32}+s_{2^21}+s_{1^5})
\\[2pt] \hline
\CM_{2,1} & q^4s_1+q^3s_1 \\
\CM_{2,2} & q^5s_2+q^4(s_2+s_{1^2})-s_2 \\
\CM_{2,3} & q^6s_3+q^5(s_3+s_{21})-q^4s_3-q(s_3+s_{21}) \\[2pt]
\hline
\end{tabular}$$
Here, if $\mu=(\mu_1\ge\dots\ge\mu_\ell)$ is a partition of $n$, $s_\mu$
denotes the associated Schur polynomial, given by the Jacobi-Trudi formula
$s_\mu=\det(h_{\mu_i+j-i})$. The above Serre polynomials are calculated in
\cite{config} (for genus $0$ and $1$) and \cite{genus2} (for genus $2$).

Applying Theorem \ref{Feynman} (and taking advantage of J. Stembridge's
symmetric function package \texttt{SF} \cite{SF} for \texttt{maple}), we
see that
\begin{align*}
\Serre(\Mbar_{3,1}) &= \Serre(\CM_{3,1}) + 3q^6 + 15q^5 +
29q^4 + 29q^3 + 16q^2 + 4q .
\end{align*}
\begin{theorem}
The Hodge polynomial $\Serre(\Mbar_{3,1})$ of $\Mbar_{3,1}$ equals
$$
q^7 + 5q^6 + 16q^5 + 29q^4 + 29q^3 + 16q^2 + 5q + 1 .
$$
\end{theorem}
\begin{proof}
We must show that $\Serre(\CM_{3,1})=q^7+2q^6+q^5+q+1$. This almost, but
not quite, follows from \cite{Looijenga}. There, it was shown that the
polynomial
$$
\sum_{i,j,k} u^iv^j t^k h^{i,j}(H^k_c(\CM_{3,1},\C)) = q^7 t^{14} + q^6
2t^{12} + q^5 t^{10} + q ( t^8 + t^7 ) + 2t^6 ,
$$
which would imply that $\Serre(\CM_{3,1})=q^7+2q^6+q^5+2$.

Unfortunately, there was a small error in that calculation; namely, two
terms were omitted from the calculation (4.8) of the Poincar\'e polynomial
of $W\backslash T'$ in case $R$ is of type $E_7$, associated to the Dynkin
subdiagrams $D_6\times A_1$ and $D_4\times A_1^2$ of $E_7$. Once they are
taken into account (the first contributes $t^7$, the second $qt^8$), we
see that
\begin{align*}
\sum_{i,j,k} u^iv^j t^k h^{i,j}(H^k(\CM_{3,1},\C)) &= q^7 t^{14} +
q^6 2t^{12} + q^5 t^{10} + q t^8 + t^7 \\ & \quad + (aqt+b)
(t^7+t^6) ,
\end{align*}
where $a$ and $b$ equal $0$ or $1$. This proves the theorem.
\end{proof}


\begin{thebibliography}{9}

\bibitem{config} E. Getzler, \emph{Resolving mixed Hodge structures of
configuration spaces.} Duke Math. J. \textbf{96} (1999), 175--203.

\bibitem{genus2} E. Getzler, \emph{Topological recursion relations in genus
$2$,} in ``Integrable systems and algebraic geometry (Kobe/Kyoto, 1997),''
World Sci. Publishing, River Edge, NJ, 1998, pp.\ 73--106.

\bibitem{modular} E. Getzler and M. Kapranov, \emph{Modular operads,}
Compositio Math. \textbf{110} (1998), 65--126.

\bibitem{Knudsen} F.F. Knudsen, \emph{The projectivity of the moduli space
of stable curves II. The stacks $\Mbar_{g,n}$,} Math. Scand. \textbf{52}
(1983), 161--189.

\bibitem{Looijenga} E. Looijenga, \emph{Cohomology of $\CM_3$ and
$\CM^1_3$,} in ``Mapping class groups and moduli spaces of Riemann surfaces
(G\"ottingen, 1991/Seattle, WA, 1991),'' Contemp. Math. \textbf{150},
Amer. Math. Soc., Providence, RI, 1993, pp.\ 205--228.

\bibitem{SF} J. Stembridge, \emph{A Maple package for symmetric functions.}
\newline \texttt{http://www.math.lsa.umich.edu/$\sim$jrs/maple.html}.

\end{thebibliography}
\end{document}